\documentclass[review]{elsarticle}
\usepackage{hyperref}
\usepackage{amsfonts}

\begin{document}

\newtheorem{thm}{Theorem}
\newtheorem{lem}[thm]{Lemma}
\newtheorem{cor}[thm]{Corollary}
\newtheorem{pro}[thm]{Proposition}
\newdefinition{dfn}[thm]{Definition}
\newproof{prf}{Proof}


\title{Excluding Degree Sequences of Cycles}
\author[cja]{Christian Joseph Altomare}
\ead{altomare@math.ohio-state.edu}
\address{The Ohio State University, 231, West 18th Avenue,
Columbus, Ohio, United States}

\begin{abstract}
In this paper, we characterize the degree sequences excluding
the degree sequence of a square
in terms of forcibly chordal graphs, and we prove several related results.
\end{abstract}

\maketitle



\begin{keyword}
degree sequence \sep graph theory \sep excluding \sep cyle \sep chordal
\sep forcibly
\end{keyword}

\section{Introduction}

\def\s{\hbox{{\rm SPLIT}}}
\def\p{\preceq}
\def\c{\circ}
\def\sq{\sqsubseteq}

Let $G$ be a finite, simple graph and let $D(G)=(d_1,\ldots,d_n)$
be its sequence of vertex degrees listed in decreasing order.
The sequence $D(G)$ is known as the {\em degree sequence} of $G$, the graph
$G$ is said to {\em realize}
$D(G)$, and $G$ is said to be a {\em realization} of $D(G)$.
We call a sequence $(d_1,\ldots,d_n)$ of nonnegative integers 
a {\em  degree sequence} if it is realized by some graph.
We make the convention that if $D$ is used without comment,
it denotes a degree sequence. Similarly for $D_1$, $D_2$, and so on.

Given $D_1$ and $D_2$, we define $D_1\preceq D_2$
to mean there is $G_1$ realizing $D_1$ and $G_2$ realizing $D_2$
such that $G_1\sqsubseteq G_2$, where $\sqsubseteq$ is the induced
subgraph relation. The reader may check
that $\preceq$ is a transitive relation on degree sequences.
We say that $G_2$ excludes $G_1$ if $G_1\not\sqsubseteq G_2$,
that $D_2$ excludes $D_1$ if $D_1\not\p D_2$,
that $D$ excludes $G$ if $D$ excludes $D(G)$,
and that $G$ excludes $D$ if $D(G)$ excludes $D$.

A set $X$ of vertices in a graph is {\em complete} ({\em anticomplete})
if every (no) two distinct vertices in $X$ are adjacent. Given disjoint
sets $X$ and $Y$, we say that $X$ is {\em complete to}
({\em anticomplete to}) $Y$ if each (no) $x$
in $X$ is adjacent to each (any) $y$ in $Y$.
A graph $S$ is {\em split} if its vertex set has a partition $(A,B)$
with $A$  complete and $B$  anticomplete.
Such partitions are called {\em split partitions}.
Split partitions of a split graph are not in general unique.
In writing a split partition $(A,B)$, we understand that $A$ is complete
and $B$ is anticomplete.

A graph is {\em chordal} if every induced cycle is a triangle.
For a graph $G$ and graph property $\cal{P}$, we say $G$ and $D(G)$
are {\em forcibly}-$\cal{P}$ if every realization of $D(G)$ has property
$\cal{P}$.
We let $G_1\coprod G_2$ denote the disjoint
union of $G_1$ and $G_2$.
We let $G[X]$  be the induced subgraph of $G$ with vertex set $X$.
We call a two edge matching $M_2$.
For other basic graph theoretic definitions and terminology, we refer the
reader to \cite{diestel_book}.

In \cite{proof_rao}, Chudnovsky and Seymour prove Rao's Conjecture;
given infinitely many degree sequences $D_1,D_2,\ldots,D_n,\ldots$,
there are positive integers $i<j$ such that $D_i\preceq D_j$.
To prove this they essentially give, for an arbitrary degree
sequence $D$, an approximate structure theorem for those
graphs excluding $D$. This general theorem is very powerful, allowing
them to resolve a nearly thirty year old conjecture. It further suggests and
leaves open a problem of independent interest: to give {\em exact} structure
theorems for {\em specific} degree sequence exclusions. That is the focus of
the current paper.
Our main result is a structure theorem characterizing degree sequences
excluding the square $C_4$.

\section{Technical Lemmas}

We recall the following folklore theorem, whose simple proof we omit.

\begin{pro} \label{split_characterization}
A graph $G$ is a split graph iff $G$ excludes~$M_2$ and all holes.
Equivalently, $G$ is split iff $G$ excludes~$M_2$, $C_4$, and $C_5$.
\end{pro}

In particular split graphs are chordal. We have the following corollary.

\begin{cor} \label{split_degree_characterization}
A degree sequence $D$ is the degree sequence of a split graph iff
$D$ excludes $M_2$, $C_4$, and $C_5$, or equivalently, iff
$D$ excludes $M_2$ and all cycles on at least $4$ vertices.
\end{cor}

By the well known characterization \cite{hammer_simeone}
of split graphs as those graphs for which some Erd\"os-Gallai
inequality \cite{erdos_gallai} is equality, we see that every realization
of a split graph is also  split; every split graph is thus forcibly split.

Let $S$ be  split with split partition $(A,B)$.
Let $H$ be an arbitrary graph.
We define $(S,A,B)\c H$ as the graph $G$ with vertex set $V(S)\cup V(H)$ formed
by joining $H$ completely to $A$ and anticompletely to $B$. This
operation is defined by R. Tyshkevish in \cite{tyshkevich}, where she states
and proves a unique decomposition theorem for finite graphs with respect
to $\c$.

\begin{lem} \label{S_or_H}
Let $n\ge 4$.
Let $S$ be a split graph with split partition $(A,B)$ and let $H$ be
an arbitrary graph. If $C_n$ is an induced subgraph
of  $(S,A,B)\c H$, then $C_n$ is an induced subgraph of $S$
or~$H$.
\end{lem}

\begin{prf}
Choose an induced $n$ point cyle $C$ in $(S,A,B)\c H$.
The set $V(C)\cap B$ is either empty or nonempty. We consider these cases.

Suppose there is $x$ in $V(C)\cap B$.
Since $d_C(x)=2$ and $x$ is only adjacent to vertices in $A$, we see
that $|V(C)\cap A|\ge 2$, so choose distinct vertices $y,z$ in $V(C)\cap A$.
If $V(C)\cap H$ contains a vertex $w$, then $\{y,z,w\}$ is the vertex set
of an induced triangle in $C$ since $A$ is both complete and complete to $H$.
By hypothesis, $C$ is a cycle on $n\ge 4$ vertices and thus has no induced
triangle. Therefore $V(C)\cap H$ is empty, showing $C$ is an induced subgraph
of $S$.

We now suppose $V(C)\cap B$ is empty. If $V(C)\cap A$ is also empty, then
$C\sqsubseteq H$, which proves the lemma. If $|V(C)\cap A|=1$, then $C$ has
at least three remaining elements, which are all contained in $H$ as
$V(C)\cap B$ is empty by hypothesis. Since $A$ is complete to $H$ it thus
follows that $d_C(x)\ge 3$, a contradiction. If $|V(C)\cap A|=2$, then
$V(C)\cap H$ must be nonempty so the two vertices of $V(C)\cap A$ together
with a vertex of $V(C)\cap H$ comprise an induced triangle in $C$, a
contradiction. If $|V(C)\cap A|\ge 3$, then any three vertices of
$V(C)\cap A$ comprise an induced triangle in $C$ since $A$ is complete,
a contradiction.

In any case, we have either that $C$ is contained in $S$, $H$, or a
contradiction, completing the proof.
\end{prf}

\begin{lem} \label{A_B_half_join}
Let~$k\ge 5$. Suppose~$D$ excludes~$C_{k-1}$, but~$D$ does
not exclude~$C_k$. Let~$G$ be a realization of~$D$ containing a
cycle~$C$ isomorphic to~$C_k$, and let~$A$ and~$B$ be the sets of
vertices of $G-C$ that are complete and anticomplete to~$C$,
respectively. Then~$G=(G[A\cup B],A,B)\c C$.
\end{lem}

\begin{prf}
We have only to show that every vertex of $G-C$ is complete or
anticomplete to~$C$, that vertices complete to~$C$ are pairwise
adjacent, and that vertices anticomplete to~$C$ are pairwise
nonadjacent.

First we show that every vertex of $G-C$ is complete or
anticomplete to~$C$. Assume not.
Then there is a vertex~$x$ outside of $C$ adjacent to some vertex $y$
of $C$ and nonadjacent to some other vertex $z$ of $C$. Let $v$
be a neighbor of $z$ in $C$ distinct from $y$. Let $K=G[C\cup x]$. Define
$K'$ as the graph obtained from $K/\{v,z\}$ by subdividing the edge
$xy$ with a new vertex $t$.
Simple checking shows that $K$ and $K'$ have the same
degree sequence. But $K'-\{x,t\}$ is isomorphic to~$C_{k-1}$.
Therefore $K'$ contains~$C_{k-1}$ as an induced subgraph.
Therefore~$K$ does not exclude $D(C_{k-1})$, and hence~$G$ does not
exclude $D(C_{k-1})$ either. This implies
that $D(C_{k-1})\le D$, contrary to hypothesis. This contradiction
shows that every vertex outside $C$ is complete or anticomplete to
$C$ as claimed.

Next, assume there are nonadjacent vertices~$x$ and~$y$, both complete
to~$C$.
Write~$C$ in cyclic order as $c_1,c_2,\ldots,c_k$.
Let $G'=G+c_1c_3-c_3x+xy-yc_1$. One may check that $D(G)=D(G')$ and that
$G'[c_1,c_3,c_4,\ldots,c_k]$ is a cycle in that
cyclic order. Therefore $G'$ contains an induced $C_{k-1}$.
We thus see that~$D$ does not exclude~$C_{k-1}$, contrary
to hypothesis. This contradiction shows~$x$ and~$y$ must be adjacent.
Since~$x$ and~$y$ are arbitrary elements of~$A$, it follows that~$A$
is complete as claimed.

Finally, 
let~$x$ and~$y$ be distinct vertices in~$B$. It is enough to
show~$x$ and~$y$ are not adjacent. Suppose they are adjacent.
Then $G[C\cup\{x,y\}]$ is isomorphic to $C_k\coprod P_2$,
which has the same degree sequence as $C_{k-1}\coprod P_3$.
Therefore~$D$ does not exclude~$C_{k-1}$, contrary to assumption.
This completes the proof.
\end{prf}

\section{The Main Results}

We now state our main theorem, from which we derive our other main results
as corollaries. A certain abuse of notation makes the statements of these
results more concise, so we make the convention that
$D=D(\s\c G)$ means that there is some split graph $S$ with split partition
$(A,B)$ such that $D$ is the degree sequence of $(S,A,B)\c G$,
and $D=D(\s)$ means $D$ is the degree sequence of a split graph.

\begin{thm} \label{gen_n_thm}
Let $n\ge 4$. A degree sequence $D$ excludes $C_n$
iff either $D=D(\s\c C_{n+1})$, $D=D(\s\c C_{n+2})$,
or  $D$ forcibly excludes each chordless cycle on at least $n$ vertices.
\end{thm}

\begin{prf}
First, if $D$ excludes all chordless cycles on $n$ or more vertices, then in
particular $D$ excludes $C_n$. If $D=D(\s\c C_{n+1})$,
then by Lemma \ref{S_or_H}, if $D(C_n)\p D$ then $C_n\sq C_{n+1}$ or
$C_n\sq S$ for some split graph $S$, a contradiction. Therefore $D$ excludes
$C_n$. Similarly if $D=D(\s\c C_{n+2})$ then $D$ excludes $C_n$.
One direction of the theorem is thus proved.

We now prove the converse. So, let $D$ exclude $C_n$.
We must show $D$ falls into one of the above three classes as claimed.

First, note $D$ excludes $C_{n+k}$ for all $k\ge 3$.
To see this, assume not. Note that $D(C_{n+k})=D(C_n\coprod C_k)$,
as $C_k$ exists since $k\ge 3$ by assumption.
Therefore $D(C_n\coprod C_k)\p D$,
so that $D$ has a realization $G$ such that $C_n\coprod C_k\sq G$.
In particular $C_n\sq G$, contrary to assumption that $D$
excludes  $C_n$. This contradiction proves our claim.

Next, we break into cases.
The first case we consider is that $D$ excludes
$C_{n+1}$ and $C_{n+2}$.
$D$ excludes $C_n$ by hypothesis,
and by the previous paragraph,
$D$ excludes $C_{n+k}$ for $k$ at least
three.
Therefore $D$ excludes all cycles on at least $n$ vertices. Therefore,
as claimed, no realization has a chordless cycle on $n$ or more vertices.

The other case is that $D$ does not exclude both $C_{n+1}$ and $C_{n+2}$.
Then $D$ has a realization $G$ containing either $C_{n+1}$ or $C_{n+2}$
as an induced subgraph.
Let $k=n+1$ if $G$ contains an induced $C_{n+1}$ and let $k=n+2$ otherwise.
In either case, $D$ excludes $C_{k-1}$ but not $C_k$.
Let $C_k=C$. Then it follows by Lemma \ref{A_B_half_join} that
$G=(G[A\cup B],A,B)\c C)$, thus completing the proof.
\end{prf}

Applying Theorem \ref{gen_n_thm} with $n=4$,
we obtain the following theorem as a corollary.

\begin{thm} \label{exclude_square}
A degree sequence $D$ excludes $C_4$ iff either
$D$ is forcibly chordal, $D=D(\s\c C_5)$, or $D=D(\s\c C_6)$.
\end{thm}

Taking complements yields the following theorem as well.

\begin{thm} \label{exclude_matching}
A degree sequence $D$ excludes $M_2$ iff
either $D$ is forcibly antichordal, $D=D(\s\c C_5)$,
or $D=D(\s\c K_{3,3})$.
\end{thm}

\begin{prf}
Just take complements, use Theorem \ref{exclude_square},
and note antichordal graphs
are the complements of chordal graphs by definition,
the pentagon is self-complementary,
and the complement of a hexagon has the same degree sequence as $K_{3,3}$.
\end{prf}

We omit proofs of the following corollaries.  Details may be found
in \cite{my_thesis}.

\begin{cor} \label{split_plus_hexagon_degree_characterization}
A degree sequence $D$ excludes $C_4$ and $C_5$ iff
$D$ is forcibly chordal or $D=D(\s\c C_6)$.
\end{cor}

\begin{cor} \label{deg_C_M}
A degree sequence $D$ excludes $M_2$ and $C_4$ iff
$D=D(\s\c C_5)$ or $D=D(\s)$.
\end{cor}

In fact, it is proved in \cite{my_thesis} that Corollary
\ref{deg_C_M} holds not only
for degree sequences, but in fact for graphs as well. More precisely, the
following proposition holds.

\begin{pro} \label{split_plus_pentagon_graph_characterization}
A graph $G$ excludes $M_2$ and $C_4$ iff $G$ is a split graph or
$G=(S,A,B)\c C_5$ for some split graph $S$ with split partition $(A,B)$.
\end{pro}

\bibliographystyle{plain}
\bibliography{mybib}

\begin{thebibliography}{1}

\bibitem{my_thesis}
Christian Altomare.
\newblock {\em Degree {S}equences, {F}orcibly {C}hordal {G}raphs and
  {C}ombinatorial {P}roof {S}ystems}.
\newblock PhD thesis, The Ohio State University, December 2009.
\newblock http://etd.ohiolink.edu.

\bibitem{proof_rao}
Maria Chudnovsky and Paul Seymour.
\newblock The proof of {R}ao's {C}onjecture on degree sequences.
\newblock In Preparation.

\bibitem{diestel_book}
Reinhard Diestel.
\newblock {\em Graph {T}heory}.
\newblock Graduate Texts in Mathematics. Springer-Verlag, New York, second
  edition, 2000.

\bibitem{erdos_gallai}
Paul Erd{\"o}s and Tibor Gallai.
\newblock {G}r{\'a}fok el{\H{o}}{\'i}rt fok{\'u} pontokkal {(}{G}raphs with
  {P}rescribed {D}egrees of {V}ertices{)}.
\newblock {\em Mat. Lapok}, 11:264--274, 1960.
\newblock In Hungarian.

\bibitem{hammer_simeone}
Peter~L. Hammer and Bruno Simeone.
\newblock The splittance of a graph.
\newblock {\em Combinatorica}, 1(3):275--284, 1981.

\bibitem{tyshkevich}
Regina Tyshkevich.
\newblock Decomposition of graphical sequences and unigraphs.
\newblock {\em Discrete Math}, 220(1-3):201--238, 2000.

\end{thebibliography}
\end{document}